\documentclass[a4paper]{elsarticle}
\usepackage{geometry}
\usepackage{amssymb}
\usepackage{amsthm}

\begin{document}
\begin{frontmatter}
\title{The vanishing ideal of a finite set of points with multiplicity structures}
\author{Na Lei, Xiaopeng Zheng, Yuxue Ren\\~~\\ {\small School of Mathematics, Jilin University}}
\date{}
\begin{abstract}

Given a finite set of arbitrarily distributed points in affine space with arbitrary multiplicity structures, we present
an algorithm to compute the reduced Gr$\ddot{\rm{o}}$bner basis of the vanishing
ideal under the lexicographic ordering. Our method discloses the
essential geometric connection between the relative position of the
points with multiplicity structures and the quotient basis of the
vanishing ideal, so we will explicitly know the set of leading terms
of elements of $I$.
We split the problem into several smaller ones which can be solved by induction
over variables and then use our new algorithm for intersection of ideals to compute the result of
the original problem. The new algorithm for intersection of ideals is mainly based on the Extended Euclidean Algorithm.
\end{abstract}
\end{frontmatter}

\section{Introduction}

To describe the problem, first we give the definitions below.

{\bfseries Definition 1:} $D\subseteq \mathbb{N}_{0}^{n}$ is called
a lower set as long as $\forall d\in D$ if $d_{i}\neq 0$, $d-e_{i}$ lies in $D$ where $e_{i}=(0, \ldots, 0, 1, 0,
\ldots, 0)$ with the 1 situated at the $i$-th position ($1\leq i\leq
n$). For a lower set $D$, we define its limiting set $E(D)$ to be
the set of all $\beta\in \mathbb{N}_{0}^{n}-D$ such that whenever
$\beta_{i}\neq 0$, then $\beta -e_{i}\in D$.

As showed in Fig.1 below, there are three lower sets and their limiting sets.
The elements of the lower sets are marked by solid circles and the elements of the
limiting sets are marked by blank circles.
\begin{center}
  \includegraphics[height=3.4cm]{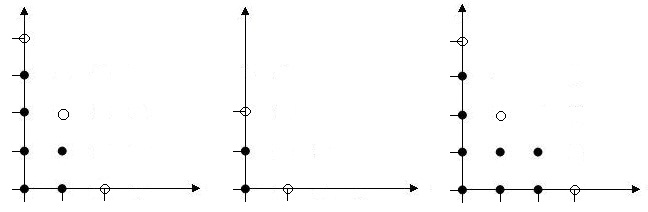}\\Fig.1: Illustration of three lower sets and their limiting sets.
\end{center}

Let $k$ be a field and $p$ be a point in the affine space
$k^{n}$, i.e. $p=(p_{1}, \ldots, p_{n})\in
k^{n}$. Let $k[X]$ be the polynomial ring over $k$, where we write
$X=(X_{1},X_{2},\ldots,X_{n})$ for brevity's sake.

{\bfseries Definition 2:} $\langle p, D\rangle$ represents a point
$p$ with multiplicity structure $D$, where $p$ is a point in affine
space $k^{n}$ and $D$ is a lower set. $\sharp D$ is called the
multiplicity of point $p$ (here we use the definition in [3]).
For each $d=(d_{1},\ldots,d_{n})\in D$, we define a
corresponding functional
$$L(f)=\frac{\partial^{d_{1}+\ldots+d_{n}}}{\partial x_{1}^{d_{1}}\ldots \partial x_{n}^{d_{n}}}f(p).$$

Hence for any given finite set of points with multiplicity structures $H=\{\langle p_{1},
D_{1}\rangle, \ldots, \langle p_{t}, D_{t}\rangle\}$, we can define $m$ functionals where $m\triangleq\sharp D_{1} +\ldots +\sharp
D_{t}$. Our aim is to
find the reduced Gr$\ddot{\rm{o}}$bner basis of the vanishing ideal $I(H)=\{f\in
k[X];L_{i}(f)=0, i=1,\ldots,m\}$ under the lexicographic ordering with $X_{1}\succ X_{2}
\succ \ldots\succ X_{n}$.

There exists an algorithm that provides a complete solution to this
problem in [4].
However, our answer for the special case of lexicographical ordering
will be in a way more transparent than the one above. %Our method obtains the quotient basis with little computation cost and get the reduced  G$\ddot{\rm{o}}$bner basis without solving any linear equations.
The ideas are summed-up as follows:

$\bullet$ Construct the reduced Gr$\ddot{\rm{o}}$bner basis of $I(H)$ and get the quotient basis $D(H)$ by induction over variables.

$\bullet$ Get the quotient basis $D(H)$ purely according to the geometric distribution of the points with multiplicity structures.

$\bullet$ Split the original problem into smaller ones which can be converted into 1 dimension lower problems and hence can be solved by induction over variables.

$\bullet$ Compute the intersection of the ideals of the smaller problems by using Extended Euclidean Algorithm.

There are several publications which have a strong connection to the work presented here. Paper [5] give a computationally efficient algorithm to get the quotient basis of the vanishing ideal over a set of points with no multiplicity structures and the authors introduce the lex game to describe the problem. Paper [6] offers a purely combinatorial algorithm to obtain the linear basis of the quotient algebra which can handle the set of points with multiplicity structures but it does not give the Gr$\ddot{\rm{o}}$bner basis. For a finite set of points with multiplicity structures,
our algorithm obtains a lower set by induction over variables and constructs the reduced Gr$\ddot{\rm{o}}$bner bases at the same time. It is only by constructing Gr$\ddot{\rm{o}}$bner basis we can prove that the lower set is the quotient basis.

One important feature of our method is the clear geometric interpretation, so in Section 2 an example together with some auxiliary pictures will be given in the first place to demonstrate this kind of feature which can make the algorithms and conclusions in this paper easier understood for us. In Section 3 and 4, some definitions and notions are given.
Section 5 and 6 are devoted to our main algorithms of computing the reduced Gr$\ddot{\rm{o}}$bner basis and the
quotient basis together with the proofs. In Section 7 we demonstrate the algorithm to compute the intersection of two ideals and some applications.

\section{Example}

First we give two different forms to represent the set of points $H$ with multiplicity structures.

For easier description, we introduce the matrix form which consists of two matrices $\langle \mathcal{P}=(p_{i,
j})_{m\times n}, \mathcal{D}=(d_{i, j})_{m\times n}\rangle$ with $\mathcal{P}_{i},
\mathcal{D}_{i}$ denoting the $i$-th row vectors of $\mathcal{P}$ and
$\mathcal{D}$ respectively. Each pair $\{\mathcal{P}_{i},
\mathcal{D}_{i}\}$ $(1\leq i\leq m)$ defines a functional in the following way.

$$L_{i}(f)=\frac{\partial^{d_{i, 1}+\ldots+d_{i, n}}}{\partial x_{1}^{d_{i, 1}}\ldots \partial x_{n}^{d_{i, n}}}f|_{x_{1}=p_{i, 1}, \ldots, x_{n}=p_{i, n}.}$$

And the functional set defined above is the same with that defined by $H$ in Section 1.

For example, given a set of three points with their
multiplicity structures $\{\langle p_{1}, D_{1}\rangle, \langle
p_{2}, D_{2}\rangle, \langle p_{3}, D_{3}\rangle\}$, where
$p_{1}=(1, 1), p_{2}=(2, 1), p_{3}=(0, 2), D_{1}=\{(0, 0), (0,
1), (1, 0)\},D_{2}=\{(0,0),(0,1),(1,0),(1,1)\}, D_{3}=\{(0, 0), (1, 0)\}$, the matrix form is like the follows.

$$
\mathcal{P}=\left(\begin{array}{cc}1&1\\1&1\\1&1\\2&1\\2&1\\2&1\\2&1\\0&2\\0&2
\end{array}\right)
,
\mathcal{D}=\left(\begin{array}{cc}0&0\\1&0\\0&1\\0&0\\1&0\\0&1\\1&1\\0&0\\1&0
\end{array}\right).
$$

For intuition's sake, we also represent the points with multiplicity structures in a more intuitive way as showed in the left picture of Fig.2 where each lower set which represents the multiplicity structure of the corresponding point $p$ is also put in the affine space with the zero element (0,0) situated at $p$. This intuitive representing form is the basis of the geometric interpretation of our algorithm.

We take the example above to show how our method works and what the geometric interpretation of our algorithm is like:

\textbf{Step 1:}
% According to the $X_{2}$ coordinates of the points, we split the original problem defined by $H=\{\langle p_{1},D_{1}\rangle, \langle p_{2}, D_{2}\rangle, \langle p_{3},D_{3}\rangle\}$ into two small ones: one problem defined by $H_{1}=\{\langle p_{1}, D_{1}\rangle, \langle p_{2}, D_{2}\rangle\}$ and the other defined by $H_{2}=\{\langle p_{3}, D_{3}\rangle\}$ which are respectively showed in the middle and the left pictures in Fig.2. In this step, we are actually to consider the $\pi$-fibres where $\pi:H\mapsto k$ such that $\langle p=(p_1,\ldots,p_n),D\rangle\in H$ is mapped to $p_n\in k$.
Define mapping $\pi:H\mapsto k$ such that $\langle p=(p_1,\ldots,p_n),D\rangle\in H$ is mapped to $p_n\in k$. So $H=\{\langle p_{1},D_{1}\rangle, \langle p_{2}, D_{2}\rangle, \langle p_{3},D_{3}\rangle\}$ consists of two $\pi$-fibres: $H_{1}=\{\langle p_{1}, D_{1}\rangle, \langle p_{2}, D_{2}\rangle\}$ and $H_{2}=\{\langle p_{3}, D_{3}\rangle\}$ as showed in the middle and the right pictures in Fig.2. Each fibre defines a new problem, so we split the original problem defined by $H$ into two small ones defined by $H_1$ and $H_2$ respectively.

\begin{center}
  \includegraphics[height=4.7cm]{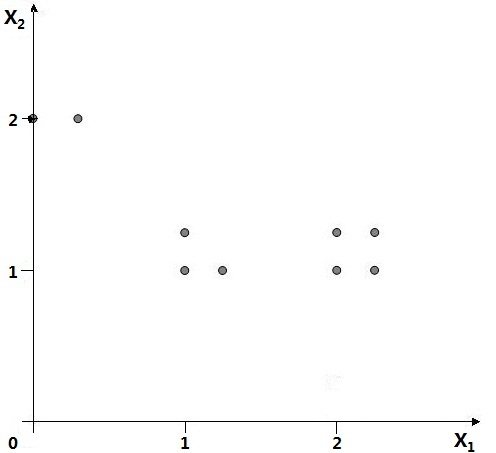}\hskip 0.1cm
  \includegraphics[height=4.7cm]{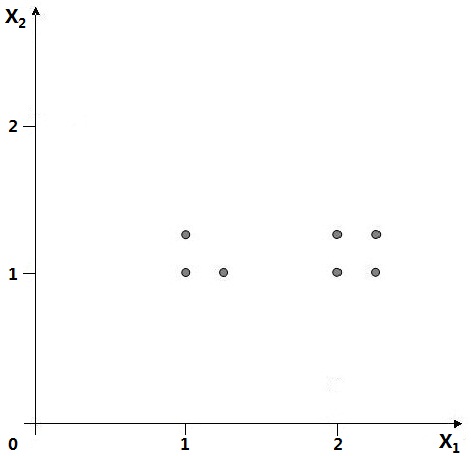}
  \includegraphics[height=4.7cm]{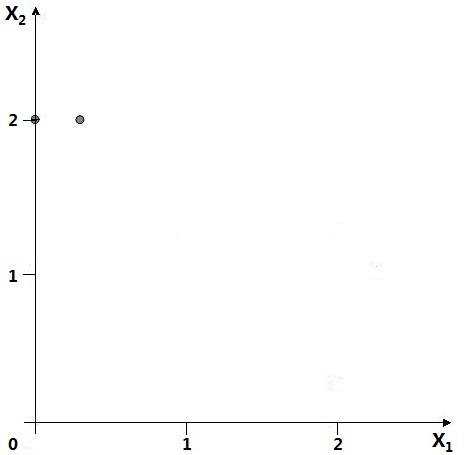}\\Fig. 2: The left picture represents $H$. The middle one is for $H_{1}$ and the right one for $H_{2}$.
\end{center}

\textbf{Step 2:} Solve the small problems. Take the problem defined
by $H_{1}$ for example.

First, it's easy to write down one element of $I(H_{1})$:
$$f_{1}=(X_{2}-1)(X_{2}-1)=(X_{2}-1)^{2}\in I(H_{1}).$$

The geometry interpretation is: we draw two lines sharing the same equation of $X_{2}-1=0$
to cover all the points as illustrated in the left picture in Fig.3 and the corresponding polynomial is $f_1$.

\begin{center}
         \includegraphics[height=4.5cm]{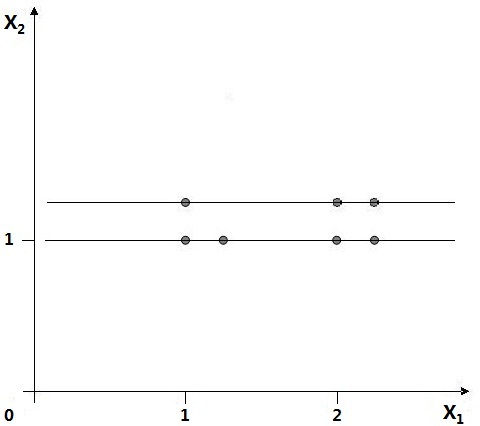}
         \includegraphics[height=4.5cm]{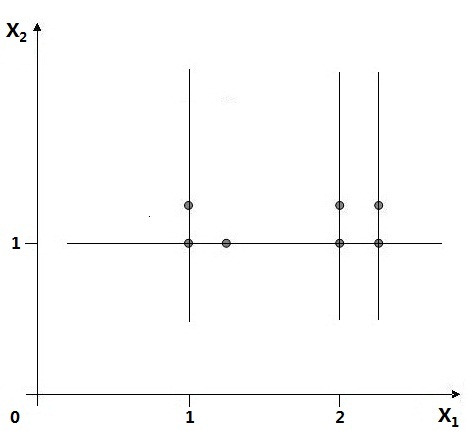}
         \includegraphics[height=4.5cm]{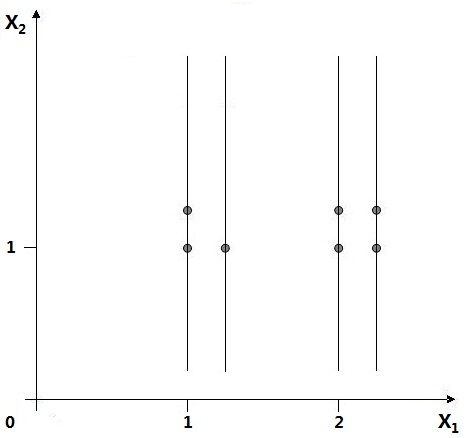}\\Fig. 3: Three ways to draw lines to cover the points.
\end{center}

According to the middle and the right pictures in Fig.2, we can
write down another two polynomials in $I(H_{1})$:
$$f_{2}=(X_{2}-1)(X_{1}-1)(X_{1}-2)^{2}~~{\rm and}~~f_{3}=(X_{1}-1)^{2}(X_{1}-2)^{2}.$$

It can be checked that $G_{1}=\{ f_{1},f_{2},f_{3}\}$ is the
reduced Gr$\ddot{\rm{o}}$bner basis of $I(H_{1})$, and the quotient basis is $\{
1,X_{1},X_{2},X_{1}X_{2},X_{1}^{2},X_{2}X_{1}^{2},X_{1}^{3}\}$. In the following,
we don't distinguish explicitly an $n$-variable monomial $X_1^{d_1}X_2^{d_2}\ldots X_n^{d_n}$ with the element $(d_1,d_2,\ldots,d_n)$ in $\mathbb{N}_{0}^{n}$, and we denote
the quotient basis of $I(H)$ by $D(H)$. Hence $D(H_{1})$ can be written
as a subset of $\mathbb{N}_{0}^{n}$:
$\{(0,0),(1,0),(0,1),(1,1),(2,0),(2,1),(3,0)\}$, i.e. a lower set, denoted
by $D^{'}$.

In fact we can get the lower set in a more direct way by pushing
the points with multiplicity structures leftward which is illustrated in the picture
below (lower set $D^{'}$ is positioned in the right part of the picture with the (0,0) element situated at point (0,1)). The elements of the lower set $D^{'}$ in the right picture in Fig.4) are marked by solid
circles. The blank circles constitute the limiting set $E(D^{'})$ and they are the leading terms of the reduced Gr$\ddot{\rm{o}}$bner basis $\{f_{1},f_{2},f_{3}\}$.

\begin{center}
    \includegraphics[height=4.7cm]{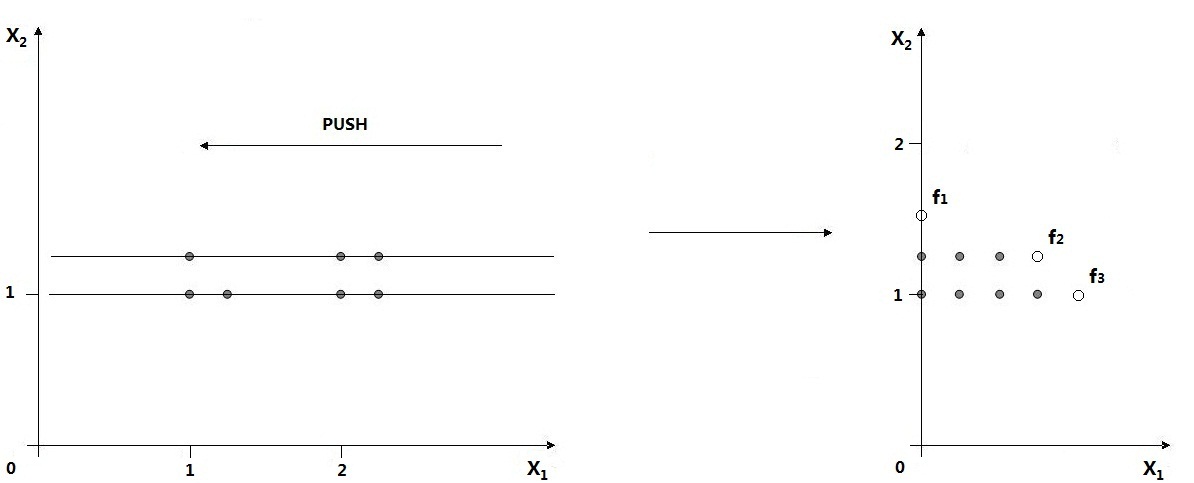}\\Fig.4: Push the points leftward to get a lower set.
\end{center}

In the same way, we can get the Gr$\ddot{\rm{o}}$bner basis $G_{2}=\{h_{1},h_{2}\}$ and a lower set $D^{''}$ for the problem defined by $H_{2}$, where $h_{1}=(X_{2}-2),h_{2}=X_{1}^{2},D^{''}=\{(0,0),(1,0)\}$.

\textbf{Step 3:} Compute the intersection of the ideals $I(H_{1})$ and $I(H_{2})$ to get the
result for the problem defined by $H$.

First, we construct a new lower set $D$ based on $D^{'},D^{''}$ in an intuitive
way: let the solid circles fall down and the elements of $D^{''}$ rest on the elements of $D^{'}$ to form a new lower set $D$ which is showed in the right part of Fig.5 and the blank circles represent the elements of the limiting set $E(D)$.

\begin{center}
    \includegraphics[height=5cm]{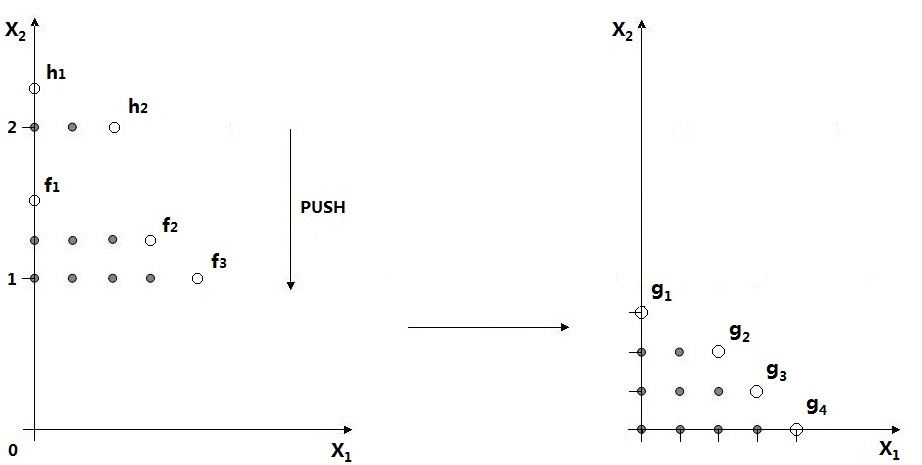}\\Fig. 5: Get the lower set $D$ based on $D^{'}$ and $D^{''}$.
\end{center}

Then we need to find $\sharp E(D)$ polynomials vanishing on $H$
with leading terms being the elements of $E(D)$. Take
$X_{1}^{3}X_{2}\in E(D)$ for example to show the general way we do it.

We need two polynomials which vanish on $H_{1}$ and
$H_{2}$ respectively, and their leading terms both have the same degrees
of $X_{1}$ with that of the desired monomial $X_{1}^{3}X_{2}$ and both have
the minimal degrees of $X_{2}$. It's easy to notice that $f_{2}$ and $X_{1}\cdot h_{2}$ satisfy the requirement and then we multiply $f_{2}$ and $X_{1}\cdot h_{2}$ with $h_{1},f_{1}$ respectively which are all univariate
polynomials of $X_{2}$ to get two polynomials $q_{1},q_{2}$
which both vanish on
$H$.
$$q_{1}=f_{2}\cdot h_{1}=(X_{2}-1)(X_{1}-1)(X_{1}-2)^{2}(X_{2}-2),$$
$$q_{2}=X_{1}\cdot h_{2}\cdot f_{1}=X_{1}^{3}(X_{2}-1)^{2}.$$

Next try to find two univariate polynomials of $X_{2}$: $r_{1},r_{2}$ such that $q_{1}\cdot
r_{1}+q_{2}\cdot r_{2}$ vanishes on $H$ (which is apparently true already) and has the desired leading
term $X_{1}^{3}X_{2}$.

To settle the leading term issue, write $q_{1},q_{2}$ as univariate polynomials of $X_{1}$.
$q_{1}=(X_{2}-2)(X_{2}-1)X_{1}^{3}-(5X_{2}^{2}-15X_{2}+10)X_{1}^{2}+(8X_{2}^{2}-24X_{2}+16)X_{1}-4X_{2}^{2}+12X_{2}-8,
q_{2}=(X_{2}-1)^{2}X_{1}^{3}$. Because $X_{2}\prec X_{1}$ and the
highest degrees of $X_{1}$ of the leading terms of $q_{1},q_{2}$ are
both $3$, we know that as long as the leading term of
$(X_{2}-2)(X_{2}-1)X_{1}^{3}\cdot r_{1}+(X_{2}-1)^{2}X_{1}^{3}\cdot
r_{2}$ is $X_{1}^{3}X_{2}$, the leading term of $q_{1}\cdot
r_{1}+q_{2}\cdot r_{2}$ is also $X_{1}^{3}X_{2}$.

$$(X_{2}-2)(X_{2}-1)X_{1}^{3}\cdot r_{1}+(X_{2}-1)^{2}X_{1}^{3}\cdot r_{2}$$
$$=X_{1}^{3}(X_{2}-1)\left((X_{2}-2)\cdot r_{1}+(X_{2}-1)\cdot r_{2}\right)$$

Obviously if and only if $(X_{2}-2)\cdot r_{1}+(X_{2}-1)\cdot r_{2}=1$ we can
keep the leading term of $q_{1}\cdot r_{1}+q_{2}\cdot r_{2}$ to be
$X_{1}^{3}X_{2}$. In this case $r_{1}=-1$ and $r_{2}=1$ will be just
perfect. In our algorithm we use Extended Euclid Algorithm to compute
$r_{1},r_{2}$.

Finally we obtain
$$g_{3}=q_{1}\cdot r_{1}+q_{2}\cdot r_{2}=(X_{2}-1)X_{1}^{3}+(5X_{2}^{2}-15X_{2}+10)X_{1}^{2}-(8X_{2}^{2}-24X_{2}+16)X_{1}+4X_{2}^{2}-12X_{2}+8$$
which vanishes on $H$ and has $X_{1}^{3}X_{2}$ as its leading term.

In the same way, we can get $g_{1}=(X_{2}-1)^{2}(X_{2}-2)$ for $X_{2}^{3}$, $g_{2}=(X_{2}-1)^{2}X_{1}^{2}$ for $X_{1}^{2}X_{2}^{2}$ and
$g_{4}=X_{1}^{4}+6(X_{2}^{2}-2X_{2})X_{1}^{3}-13(X_{2}^{2}-2X_{2})X_{1}^{2}+12(X_{2}^{2}-2X_{2})X_{1}-4(X_{2}^{2}-2X_{2})$ for $X_{1}^{4}$. In fact we need to compute $g_{1},~g_{2},~g_{3}$ and $g_{4}$ in
turn according to the lexicographic order because we need reduce
$g_{2}$ by $g_{1}$, reduce $g_{3}$ by $g_{2}$ and $g_{1}$, and reduce $g_4$ by $g_1$, $g_2$ and $g_3$.

The reduced polynomial set can be proved in Section 6 to be
the reduced Gr$\ddot{\rm{o}}$bner basis of the intersection of two ideals which is exactly the vanishing ideal over $H$, and $D$ is the
quotient basis.

\section{Notions}

First, we define the following mappings.

$~~~~proj:\mathbb{N}_{0}^{n} \longrightarrow k$

$~~~~~~~~~~~(d_{1},\ldots,d_{n})\longrightarrow d_{n}$.

$~~~~\widehat{proj}:\mathbb{N}_{0}^{n}\longrightarrow \mathbb{N}_{0}^{n-1}$

$~~~~~~~~~~~(d_{1},\ldots,d_{n})\longrightarrow (d_{1},\ldots,d_{n-1})$.

$~~~~embed_{c}:\mathbb{N}_{0}^{n-1}\longrightarrow \mathbb{N}_{0}^{n}$

$~~~~~~~~~~~(d_{1},\ldots,d_{n-1})\longrightarrow (d_{1},\ldots,d_{n-1},c)$.

Let $D\subset N_0^{n}$, and naturally we define $\widehat{proj}(D)=\{\widehat{proj}(d)|d\in D\}$, and $embed_{c}(D^{'})=\{embed_{c}(d)|d\in D^{'}\}$ where $D^{'}\subset N_0^{n-1}$. In fact we can apply these mappings to any set $O\subset k^{n}$ or any matrix of $n$ columns, because there is no danger of confusion. For example, let $M$ be a matrix of $n$ columns, and $\widehat{proj}(M)$ is a matrix of $n-1$ columns with the first $n-1$ columns of $M$ reserved and the last one eliminated.

The $embed_c$ mapping embeds an $n-1$ dimensional lower set into the $n$ dimensional space.
When the $embed_c$ operation parameter $c$ is zero, we can get an $n$ dimensional lower set by
mapping each element $d=(d_{1},\ldots,d_{n-1})$ to $d=(d_{1},\ldots,d_{n-1},0)$ as showed below.

\begin{center}
    \includegraphics[height=4.3cm]{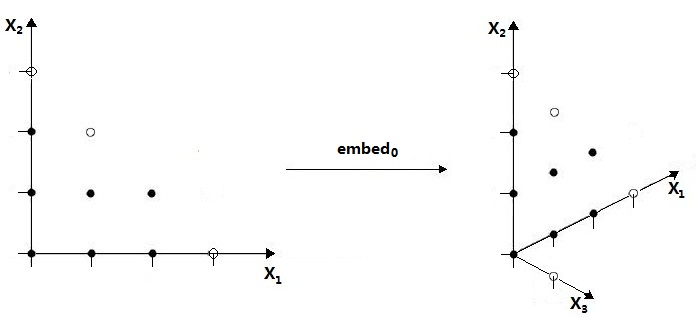}\\Fig. 6: Embed the lower set in 2-D space into 3-D space with parameter $c=0$.
\end{center}

Blank circles represent the elements of the limiting sets. Note that after the $embed_c$ mapping,
there is one more blank circle. In this case, the limiting set is always increased by one element $(0,\ldots,0,1)$.

In the case the $embed_c$ operation parameter $c$ is not zero, it is obvious that what we got is not a lower set
any more. But there is another intuitive fact we should realize.

\textbf{Theorem 1:} $D_{0},D_{1},\ldots,D_{k}$ are $n-1$ dimensional lower sets, and
$D_{0}\supseteq D_{1}\supseteq \ldots \supseteq D_{k}$. Let $\hat{D}_{i}=embed_{i}(D_{i}),i=0,\ldots,k.$
Then $D=\bigcup_{i=0}^{k}\hat{D}_{i}$ is an $n$ dimensional lower set, and $E(D)\subseteq C$ where
$C=\bigcup_{i=0}^{k}embed_{i}(E(D_{i}))\bigcup \{(0,\ldots,0,k+1)\}$.

\textbf{Proof:}
First to prove $D$ is a lower set. $\forall d\in D,$ let $i=proj(d)$, then $d\in \hat{D}_{i}$ i.e. $\widehat{proj}(d)\in\widehat{proj}(\hat{D}_i)=D_i$. Because $D_{i}$ is a lower set, hence for $j=1,\ldots,n-1,$ if $d_{j}\neq 0$, then $\widehat{proj}(d)-\widehat{proj}(e_{j})\in D_{i}$ where $e_{j}=(0, \ldots, 0, 1, 0,
\ldots, 0)$ with the 1 situated at the $j$-th position. So $d-e_{j}\in \hat{D}_{i}\subseteq D$. For $j=n$, if $i=0$, then we are finished. Else there must be $d-e_{n}\in\hat{D}_{i-1}\subseteq D$. Because if $d-e_{n}\notin\hat{D}_{i-1}$, we have $\widehat{proj}(d)\notin D_{i-1}$. Since we already have $\widehat{proj}(d)\in D_{i}$, this is contradictory to $D_{i}\subseteq D_{i-1}$.

Second, $\forall d\in E(D)$, $\widehat{proj}(d)\notin D_{i},i=0,\ldots,k$. If $\widehat{proj}(d)$ is a zero tuple, then $d_{n}$ must be $k+1$, that is $d\in C.$ Else we know $d_{n}<k+1$. If $d_{j}\neq 0,~j=1,\ldots,n-1$
, then $d-e_{j}\in embed_{d_{n}}(D_{d_{n}})$. Then
$\widehat{proj}(d)-\widehat{proj}(e_{j})\in D_{d_{n}}$, that is $\widehat{proj}(d)\in E(D_{d_{n}})$. Finally with the $embed_{d_{n}}$ operation we have $d\in embed_{d_{n}}(E(D_{d_{n}}))$ where $d_{n}<k+1$. So $d\in C$.

\section{Addition of lower sets}

In this section, we define the addition of lower sets which is the same with that in [2], the following paragraph and Fig.7 are basically excerpted from that paper with a little modification of expression.

To get a visual impression of what the addition of lower sets dose, look at the example in Fig.7. What is depicted there can generalizes to arbitrary lower sets $D_{1}$ and $D_{2}$ in arbitrary dimension $n$ and can be described as follows. Draw a coordinate system of $\mathbb{N}_{0}^{n}$ and insert $D_{1}$. Place a translate of $D_{2}$ somewhere on the $X_2$-axis. The translate has to be sufficiently far out, so that $D_{1}$ and the translate $D_{2}$ do not intersect. Then take the elements of the translate of $D_{2}$ and drop them down along the $X_2$-axis until they lie on top of the elements of $D_{1}$. The resulting lower set is denoted by $D_{1}+D_{2}$.

\begin{center}
    \includegraphics[height=4.7cm]{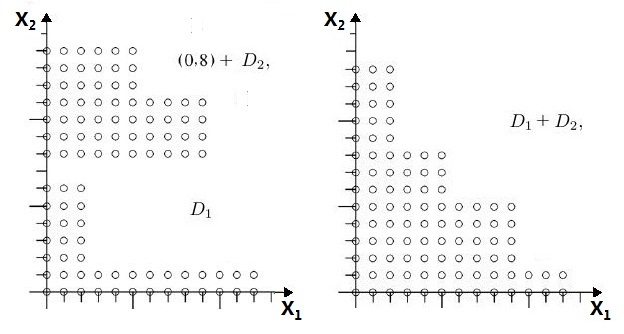}\\Fig. 7: Addition of $D_{1}$ and $D_{2}$.
\end{center}

Intuitively, we define algorithm \textbf{AOL} to realize the addition of lower sets.

Algorithm \textbf{AOL:} Given two $n$ dimensional lower sets
$D_{1},D_{2}$, determine another lower set as the addition of
$D_{1},D_{2}$, denoted by \textbf{$D:=D_{1}+D_{2}$}.

[step 1]: $D:=D_{1}$;

[step 2]: If $\sharp D_{2}=0$ return $D$. Else pick $a\in
D_{2},D_{2}:=D_{2}\setminus\{a\}.$

[step 2.1]: If $\sharp (D\bigcup \{a\})$=$\sharp D$, add the last
coordinate of $a$ with $1$. Go to [step 2.1]. Else
$D:=D\bigcup\{a\}$, go to [step 2].

Given $n$ dimensional lower sets $D_{1},D_{2},D_{3}$, the addition we
defined satisfies:

$(a)~D_{1}+D_{2}=D_{2}+D_{1},$

$(b)~(D_{1}+D_{2})+D_{3}=D_{1}+(D_{2}+D_{3}),$

$(c)~D_{1}+D_{2}$ is a lower set,

$(d)~\sharp(D_{1}+D_{2})=\sharp D_{1}+\sharp D_{2}.$

These are all the same with that in [2]. And the proof can be referred to it.

As implied in the example of Section 2, when we want to get a polynomial with leading term $d_{3}$ showed in the right part of Fig.8, we need two polynomials with the leading terms $d_{1},d_{2}$ which are not the elements of the lower sets and have the same degrees of $X_{1}$ as $d_3$ and the minimal degrees of $X_{2}$ as showed in the left part of Fig.8. In other words, $d_{1}\notin D_{1},~d_{2}\notin D_{2},~\widehat{proj}(d_{1})=\widehat{proj}(d_{2})=\widehat{proj}(d_{3})$, $proj(d_{1})+proj(d_{2})=proj(d_{3})$. It's easy to understand that these equations hold for the addition of three or even more lower sets.

\begin{center}
    \includegraphics[height=5cm]{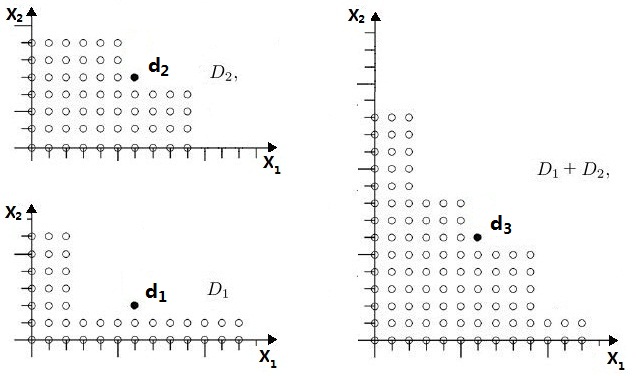}\\Fig.8: $\widehat{proj}(d_{1})=\widehat{proj}(d_{2})=\widehat{proj}(d_{3}),~proj(d_{1})+proj(d_{2})=proj(d_{3})$.
\end{center}

We use algorithm \textbf{GLT} to get the leading terms $d_{1}$ and $d_{2}$ from $d_3$ respectively.

Algorithm \textbf{GLT:} Given $a\in\mathbb{N}_{0}^{n}$, and an $n$ dimensional lower set $D$ satisfying $a\notin D$. Determine
another $r=(r_1,\ldots,r_n)\in \mathbb{N}_{0}^{n}$ which satisfies that $r\notin D$, $\widehat{proj}(r)=\widehat{proj}(a)$ and $(r_1,\ldots,r_{n-1},r_n-1)\in D$, denoted by \textbf{$r:=GLT(a,D)$}.

[step 1]: Initialize $r$ such as $\widehat{proj}(r)=\widehat{proj}(a)$ and $proj(r)=0$.

[step 2]: if $r\notin D$, return r, else $r_n:=r_n+1$, go to [step 2].

Then $d_{1}=GLT(d_{3},D_{1}),~d_{2}=GLT(d_{3},D_{2}).$

\textbf{Definition 3:} For any $f\in k[X]$, view it as an element in $k(X_{n})[X_{1},\ldots,X_{n-1}]$ and define $LC_n(f)$ to be the leading coefficient of $f$ which is an univariate polynomial of $X_{n}$.

Algorithm \textbf{GLP:} $D$ is an $n$ dimensional lower set, $a\in \mathbb{N}_{0}^{n}$ and $a\notin D$, $G:=\{f\in k[X];\exists~ ed\in E(D),s.t.$ the leading term of $f$ is $ed~\}$, algorithm \textbf{GLP} returns a polynomial $p$ in the ideal $\langle G\rangle$ whose leading term is $GLT(a,D).$ Denoted by $p:=GLP(a,D,G).$

[step 1:] $c:=GLT(a,D)$.

[step 2:] Select $c^{'}\in E(D),s.t.~c^{'}$ is a factor of $c.~~d:=\frac{c}{c^{'}}$.

[step 3:] $p:=f_{c^{'}}\cdot d$ where $f_{c^{'}}$ is an element of $G$ whose leading term is $c^{'}$.

\textbf{Remark 1:} $LC_n(f_{c^{'}})=LC_n(p)$ in [step 3]. Since $c$ has the minimal degree of $X_n$ according to algorithm \textbf{GLT}, there exists no element $c^{''}\in E(D)$ which is a factor of $c$ satisfying $proj(c^{''})<proj(c)$. Hence monomial $d$ in the algorithm does not conclude the variable $X_n$.

\section{Associate a lower set $D(H)$ to a set of points $H$ with multiplicity structures}
For any given set of $n$ dimensional points $H$ with multiplicity structures, we can construct an $n$ dimensional lower set $D(H)$ by induction.

\textbf{Univariate case:}
$H=\{\langle p_{1},D_{1}\rangle,\ldots,\langle p_{t},D_{t}\rangle\}$, then the lower set is $D(H)=\{0,1,\ldots,\sum_{i=1}^{t}\sharp D_{i}\}$.

To pass from $n-1$ to $n$ ($n\geq 2$), we first solve a \textbf{Special case}.

\textbf{Special case:}
$H=\{\langle p_{1},D_{1}\rangle,\ldots,\langle p_{t},D_{t}\rangle\}$ is a set of $n$ dimensional points with multiplicity structures where all the points share the same $X_{n}$ coordinates. Write $H$ in matrix form as $\langle \mathcal{P},\mathcal{D}\rangle$ and all the entries in the last column of matrix $\mathcal{P}$ have the same values. Classify the row vectors of $\langle\mathcal{P},\mathcal{D}\rangle$ to get $\{\langle \mathcal{P}_{0},\mathcal{D}_{0}\rangle,\ldots,\langle \mathcal{P}_{w},\mathcal{D}_{w}\rangle\}$ according to the values of the entries in the last column of matrix $\mathcal{D}$ and we guarantee the corresponding relationship between the row vectors of matrix $\mathcal{P}$ and matrix $\mathcal{D}$ holds in $\langle \mathcal{P}_{i},\mathcal{D}_{i}\rangle$ ($0\leq i\leq w$). All the entries in the last column of $\mathcal{D}_{i}$ are the same $i$ and the entries of the last column of $\mathcal{P}_{i}$ stay the same too. Then eliminate the last columns of $\mathcal{P}_{i}$ and $\mathcal{D}_{i}$ to get $\langle\widehat{proj}(\mathcal{P}_{i}),\widehat{proj}(\mathcal{D}_{i})\rangle$ which represents a set of $n-1$ dimensional points with multiplicity structures, by induction we get a lower set $\hat{D}_{i}$ in $n-1$ dimensional space. Then we set
$$D(H)=\bigcup_{i=0}^{w}embed_{i}(\hat{D}_{i}).$$

Next we deal with the \textbf{General case}.

\textbf{General case:}
$H=\{\langle p_{1},D_{1}\rangle,\ldots,\langle p_{t},D_{t}\rangle\}$ is a set of $n$ dimensional points with multiplicity structures. Split the set of points: $H=H_{1}\bigcup H_{2}\bigcup\ldots\bigcup H_{s}$. The points of $H_{i}$ are in the same $\pi$-fibre, i.e. they have the same $X_{n}$ coordinates $c_{i}$, $i=1,\ldots,s$,and $c_{i}\neq c_{j},\forall i,j=1,\ldots,s,i\neq j.$ According to the \textbf{Special case}, for each $i=1,\ldots,s$, we can get a lower set $D(H_{i})$, then we set
$$D(H)=\sum_{i=1}^{s}D(H_{i}).$$

We now proof $D(H)$ is a lower set although it is easy to understand as long as the geometric interpretation involves. Since it is obviously true for \textbf{Univariate case}, induction over dimension would be helpful for the proof.

\textbf{Proof:} Assume $D(H)$ is a lower set for the $n-1$ dimensional situation and now we prove the conclusion for $n$ dimensional situation ($n\geq 2$).

First to prove $D(H)$ of the \textbf{Special case} is a lower set.

We claim that $\langle\widehat{proj}(\mathcal{P}_{i}),\widehat{proj}(\mathcal{D}_{i})\rangle$
represents an $n-1$ dimensional set of points with multiplicity structures ($i=0,\ldots,w$). For any $D\subset \mathbb{N}_{0}^{n}$, define $F_{a}(D)=\{d\in D|~proj(d)=a\}.$ Let $U=\{u|u\in\{1,\ldots,t\},F_i(D_u)\neq \varnothing\}.$ %$\bigcup_{u\in U}\widehat{proj}(p_{u})$ will be the set of row vectors of $\widehat{proj}(\mathcal{P}_{i})$, and $\bigcup_{u\in U}\widehat{proj}(F_{i}(D_{u}))$ will be the set of row vectors in $\widehat{proj}(\mathcal{D}_{i})$.
 So $\langle\widehat{proj}(\mathcal{P}_{i}),\widehat{proj}(\mathcal{D}_{i})\rangle$ can be written in the form of $\{\langle \widehat{proj}(p_{u}), \widehat{proj}(F_{i}(D_{u})) \rangle | u\in U\}$. Apparently $\widehat{proj}(F_{i}(D_{u}))$ is an $n-1$ dimensional lower set and can be viewed as the multiplicity structure of the point $\widehat{proj}(p_{u})$. Hence $\langle\widehat{proj}(\mathcal{P}_{i}),\widehat{proj}(\mathcal{D}_{i})\rangle$
is an $n-1$ dimensional set of points with multiplicity structures.

What's else, we assert $\widehat{proj}(\mathcal{P}_{j})$ is a sub-matrix of $\widehat{proj}(\mathcal{P}_{i}),$ and $\widehat{proj}(\mathcal{D}_{j})$ is a sub-matrix of $\widehat{proj}(\mathcal{D}_{i}),0\leq i<j\leq w.$ Because of the corresponding relationship between the row vectors in $\mathcal{P}$ and $\mathcal{D}$, we need only to prove $\widehat{proj}(\mathcal{D}_{j})$ is a sub-matrix of $\widehat{proj}(\mathcal{D}_{i})$. If it is not true, there exists a row vector $g$ of $\widehat{proj}(\mathcal{D}_{j})$ which is not a row vector of $\widehat{proj}(\mathcal{D}_{i})$. That is, there exists $b$ ($1\leq b\leq t$) such that $embed_j(g)$ is an element of the lower set $D_{b}$, and $embed_i(g)$ is not included in any lower set $D_{a}$ ($1\leq a\leq t$). However since $i<j$ and $embed_j(g)\in D_b$, $embed_i(g)$ must be included in $D_b$.
Hence our assertion is true.

Since $\widehat{proj}(\mathcal{P}_{j})$ is a sub-matrix of $\widehat{proj}(\mathcal{P}_{i}),$ and $\widehat{proj}(\mathcal{D}_{j})$ is a sub-matrix of $\widehat{proj}(\mathcal{D}_{i}),0\leq i<j\leq w.$ According to the assumption of induction and the way we construct $D(H)$, we have $\hat{D}_{i}\supseteq \hat{D}_{j},0\leq i<j\leq w,$ where $\hat{D}_{i},\hat{D}_{j}$ are both lower sets. Based on the \textbf{Theorem 1} in Section 3, $D(H)=\bigcup_{i=0}^{w}embed_{i}(\hat{D}_{i})$ is a lower set, and $E(D(H))\subseteq \bigcup_{i=0}^{w} embed_{i}(E(\hat{D}_{i}))\bigcup \{(0,\ldots,0,w+1)\}$.

Then to prove $D(H)$ of \textbf{General case} is a lower set. Since $D(H_{i}),i=1,\ldots,s$ are lower sets, and the addition of lower sets is also a lower set according to Section 4, $D(H)$ is obviously a lower set. The proof is finished.

\section{Associate a set of polynomials $poly(H)$ to $D(H)$}
For every lower set constructed during the induction procedure showed in the last section, we associate a set of polynomials to it.

We begin with the univariate case as we did in the last section.

\textbf{P-univariate case:}
$H=\{\langle p_{1},D_{1}\rangle,\ldots,\langle p_{t},D_{t}\rangle\}$, and $D(H)=\{0,1,\ldots,\sum_{i=1}^{t}\sharp D_{i}\}$. The set of polynomials associated to $D(H)$ is $poly(H)=\{\prod_{i=1}^{t}(X_{1}-p_{i})^{\sharp D_{i}}\}$.

Apparently, $poly(H)$ of \textbf{P-univariate case} satisfies the following \textbf{Assumption}.

\textbf{Assumption:} For any given $n-1~(n>1)$ dimensional set of points $H$ with multiplicity structures, there are the following conclusions. For any $\lambda\in E(D(H))$, there exists a polynomial $f_{\lambda}\in k[X]$ where $X=(X_{1},\ldots,X_{n-1})$ such that

$\bullet$ The leading term of $f_{\lambda}$ under lexicographic ordering is $X^{\lambda}$.

$\bullet$ The exponents of all lower terms of $f_{\lambda}$ lies in $D(H)$.

$\bullet$ $f_{\lambda}$ vanishes on $H$.

$\bullet$ $poly(H)=\{f_{\lambda}|\lambda\in E(D(H))\}$.

When we construct the set of polynomials $poly(H)$, we should make sure the assumption always holds. Now let us consider the $n~(n>1)$ dimensional situation and still begin with the special case.

\textbf{P-Special case:} Given a set of points with multiplicity structures
$H=\{\langle p_{1},D_{1}\rangle,\ldots,\langle p_{t},D_{t}\rangle\}$
or in matrix form $\langle\mathcal{P}=(p_{ij})_{m\times
n},\mathcal{D}=(d_{ij})_{m\times n}\rangle$. All the given points have the same $X_{n}$ coordinates, i.e.
the entries in the last column of $\mathcal{P}$ are the same. We compute $poly(H)$ following the steps below.

[step 1]: $c:=p_{1n};$ $w=max\{d_{in};i=1,\ldots,m\}.$

[step 2]: $\forall i=0,\ldots,w$, define $\mathcal{SD}_{i}$ as a
sub-matrix of $\mathcal{D}$ containing all the row vectors whose last
coordinates equal to $i$. Extract the corresponding row vectors of $\mathcal{P}$
to form matrix $\mathcal{SP}_{i}$, and the corresponding relationship between the
row vectors in $\mathcal{P}$ and $\mathcal{D}$ holds for $\mathcal{SP}_{i}$ and $\mathcal{SD}_{i}$.

[step 3]: $\forall i=0,\ldots,w$, eliminate the last columns of $\mathcal{SP}_{i}$ and $\mathcal{SD}_{i}$ to get $\langle\tilde{\mathcal{SP}_{i}},\tilde{\mathcal{SD}_{i}}\rangle$ which represents a set of points in $n-1$ dimensional space with multiplicity structures.
According to the induction assumption, we have the polynomial set $\tilde{G}_{i}=poly(\langle \tilde{\mathcal{SP}_{i}}$,$\tilde{\mathcal{SD}_{i}} \rangle)$ associated to the lower set $\tilde{D}_{i}=D(\langle \tilde{\mathcal{SP}_{i}}$,$\tilde{\mathcal{SD}_{i}} \rangle)$.

[step 4]: $ D:=\bigcup_{i=0}^{w}embed_{i}(\tilde{D}_{i}).$ Multiply every element of
$\tilde{G}_{i}$ with $(X_{n}-c)^{i}$ to get $G_{i}.$
$\tilde{G}:=\bigcup_{i=0}^{w}G_{i}\bigcup \{(X_{n}-c)^{w+1}\}$.

[step 5]: Eliminate the polynomials in $G$ whose leading term is not
included in $E(D)$ to get $poly(H)$.

\textbf{Theorem 2:} The $poly(H)$ got in \textbf{P-Special case} satisfies the \textbf{Assumption}.

\textbf{Proof:} According to the Section 5, $\langle \tilde{\mathcal{SP}_{i}}$,$\tilde{\mathcal{SD}_{i}} \rangle$ represents an $n-1$ dimensional set of points with multiplicity structures for $i=0,\ldots,w.$ And $\tilde{D}_{j}\supseteq \tilde{D}_{i},0\leq j\leq i\leq w$. $D$ is a lower set and $E(D)\subseteq \bigcup_{i=0}^{w}embed_{i}(E(\tilde{D}_{i}))\bigcup \{(0,\ldots,0,w+1)\}$.

For $\lambda=(0,\ldots,0,w+1)\in E(D)$, we have $f_\lambda=(X_{n}-c)^{w+1}$. It is easy to check that it satisfies the first three terms of the \textbf{Assumption}.

For any other element $ed$ of $E(D)$, $\exists k~s.t.~ed\in embed_{k}E(\tilde{D}_{k})$. So let $\tilde{ed}$ be the element in $E(\tilde{D}_{k})$ such that $ed=embed_{k}(\tilde{ed})$. We have $f_{\tilde{ed}}$ vanishes on $\langle\tilde{\mathcal{SP}_{k}},\tilde{\mathcal{SD}_{k}}\rangle$ whose leading term is $\tilde{ed}\in E(\tilde{D}_{k})$ and the lower terms belong to $\tilde{D}_{k}$.
According to the algorithm $f_{ed}=(X_{n}-c)^{k}\cdot f_{\tilde{ed}}\in poly(H)$ .

First it is easy to check that the leading term of $f_{ed}$ is $ed$ since $ed=embed_{k}(\tilde{ed})$.

Second, the lower terms of $f_{ed}$ are all in the set $S=\bigcup_{j=0}^{k}embed_{j}(\tilde{D}_{k})$ because all the lower terms of $f_{\tilde{ed}}$ are in the set $\tilde{D}_{k}$. $\tilde{D}_{0}\supseteq \tilde{D}_{1}\supseteq \ldots \tilde{D}_{k}$, so $embed_j(\tilde{D}_{k})\subset embed_j(\tilde{D}_{j})~(0\leq j\leq k)$, hence $S\subseteq D=\bigcup_{j=0}^{w}embed_{j}(\tilde{D}_{j})$ and the second term of the \textbf{Assumption} is satisfied.

Third, we are going to prove that $f_{ed}$ vanishes on all the functionals defined by $\langle \mathcal{P},\mathcal{D}\rangle$, i.e. all the functionals defined by $\langle \mathcal{SP}_i,\mathcal{SD}_i\rangle~(i=0,\ldots,w).$

When $i\neq k$, we write all the functionals defined by $\langle \mathcal{SP}_i,\mathcal{SD}_i\rangle$ in this form: $L^{'}\cdot \frac{\partial^{i}}{\partial X_n^{i}}|_{X_n=c}$ where $L^{'}$ is an $n-1$ variable functional. Since $f_{ed}=(X_{n}-c)^{k}\cdot f_{\tilde{ed}}$, apparently $f_{ed}$ vanishes on these functionals.

For $i=k$, denote by $L$ the functionals defined by $\langle\tilde{\mathcal{SP}_{i}},~\tilde{\mathcal{SD}_{i}}\rangle$, and $f_{\tilde{ed}}$ vanishes on $L$. All the functionals defined by $\langle \mathcal{SP}_k,\mathcal{SD}_k\rangle$ can be written in this form: $L^{''}\cdot \frac{\partial^{k}}{\partial X_n^{k}}|_{X_n=c}$ where $L^{''}\in L$. Since  $f_{ed}=(X_{n}-c)^{k}\cdot f_{\tilde{ed}}$, apparently $f_{ed}$ vanishes on these functionals.

So $f_{ed}$ vanish on $H$, and $f_{ed}$ satisfies the first three terms of the \textbf{Assumption}.

In summary $poly(H)$ satisfies the \textbf{Assumption}, and we finish the proof.

\textbf{Remark 2:} For $f_{\lambda}\in poly(H),\lambda\in E(D)$ where $poly(H)$ is the result got in the algorithm above, we have the conclusion that $LC_n(f_{\lambda})=(X_{n}-c)^{proj(\lambda)}$.

\textbf{P-General case:}
Given a set of points with multiplicity structures
$H$ or in matrix form $\langle\mathcal{P}=(p_{ij})_{m\times
n},\mathcal{D}=(d_{ij})_{m\times n}\rangle$, we are going to get $poly(H)$.

[step 1]: Write $H$ as $H=H_{1}\bigcup H_{2}\bigcup\ldots\bigcup H_{s}$ where $H_{i}~(1\leq i\leq s)$ is a $\pi$-fibre ($\pi:H\mapsto k$ such that $\langle p=(p_1,\ldots,p_n),D\rangle\in H$ is mapped to $p_n\in k$) i.e. the points of $H_{i}$ have the same $X_{n}$ coordinates $c_{i}$, $i=1,\ldots,s$,and $c_{i}\neq c_{j},\forall i,j=1,\ldots,s,i\neq j.$

[step 2]: According to the \textbf{P-Special case}, we have $D^{'}_{i}=D(H_i), G_{i}=poly(H_{i})$. Write $H_{i}$ as $\langle \mathcal{P}_{i},\mathcal{D}_{i}\rangle$,
and define $w_{i}$ as the maximum value of the elements in the last column of $\mathcal{D}_{i}$.

[step 3]: $D:=D^{'}_{1}, G:=G_{1},i:=2$.

[step 4]: If $i>s$, go to [step 5]. Else

[step 4.1]: $D:=D+D^{'}_{i};$ $\hat{G}:=\varnothing$. View $E(D)$ as a monomial set $MS:=E(D)$.

[step 4.2]: If $\sharp MS= 0$, go to [step 4.7], else select the minimal element of $MS$ under lexicographic ordering, denoted by $LT$. $
MS:=MS\setminus\{LT\}$.

[step 4.3]:

$$f_{1}:=GLP(LT,D,G),f_{2}:=GLP(LT,D_{i}^{'},G_{i}).$$

~~~~~~~~~~~~~~~~~~~~~~~~~$v_{k}:=proj(g_{k})$, where $g_{k}:=GLT(LT,D_{k}^{'})$, $k=1,\ldots,i$.

[step 4.4]:
$$q_{1}:=f_{1}\cdot (X_{n}-c_{i})^{w_{i}+1};~~q_{2}:=f_{2}\cdot \prod_{k=1}^{i-1}(X_{n}-c_{k})^{w_{k}+1}.$$

$$pp1:=(X_{n}-c_{i})^{w_{i}+1-v_{i}};~~pp2:=\prod_{k=1}^{i-1}(X_{n}-c_{k})^{w_{k}+1-v_{k}}.$$

[step 4.5]: Use Extended Euclidean Algorithm to compute $r_{1}$ and
$r_{2}$ s.t. $r_{1}\cdot pp_{1}+r_{2}\cdot pp_{2}=1$.

[step 4.6]: $f:=r_{1}\cdot q_{1}+r_{2}\cdot q_{2}$. Reduce $f$ with the
elements in $\hat{G}$ to get $f^{'}$; $\hat{G}:=\hat{G}\bigcup\{f^{'}\}.$ Go to [step
4.2].

[step 4.7]: $G:=\hat{G}.$ $i:=i+1.$ Go to [step 4].

[step 5]: $poly(H):=G$.

\textbf{Theorem 3:} The $poly(H)$ got in \textbf{p-General case} satisfies the \textbf{Assumption}.

\textbf{proof:} We need only to prove the situation that $s\geq 2$ in [step 1].

For $i=2$, $D=D_{1}^{'}+D_{2}^{'}. ~~\forall ed\in E(D)$, $v:=proj(ed)$ and $X_0:=\frac{X^{ed}}{X_n^{v}}$. According to Section 4, we have $v=v_{1}+v_{2}$. Based on the \textbf{Remark 1} and \textbf{Remark 2}, $f_{1}$ and $f_{2}$ can be written as polynomials of $k(X_{n})[X_{1},\ldots,X_{n-1}]:$ $f_{1}=X_{0}\cdot (X_{n}-c_{1})^{v_{1}}+the~rest$ and $f_{2}=X_{0}\cdot (X_{n}-c_{2})^{v_{2}}+the~rest$ and none of the monomials in $the~rest$ is greater than or equal to $X_{0}.$
Because $f_{1}$ and $(X_{n}-c_{1})^{w_{1}+1}$ vanish on $H_{1}$, $f_{2}$ and $(X_{n}-c_{2})^{w_{2}+1}$ vanish on $H_{2}$, we know that $q_{1}=f_{1}\cdot (X_{n}-c_{2})^{w_{2}+1}$ and $q_{2}=f_{2}\cdot (X_{n}-c_{1})^{w_{1}+1}$ both vanish on $H_{1}\bigcup H_{2}$. Then $f$ vanishes on $H_{1}\bigcup H_{2}$ where $f=r_{1}\cdot q_{1}+r_{2}\cdot q_{2}$.

$\qquad~~f=X_{0}\cdot (X_{n}-c_{1})^{v_{1}}\cdot (X_{n}-c_{2})^{v_{2}}(r_{1}\cdot (X_{n}-c_{2})^{w_{2}+1-v_{2}}+r_{2}\cdot (X_{n}-c_{1})^{w_{1}+1-v_{1}})+the~ rest$

$\qquad~~~~=X_{0}\cdot (X_{n}-c_{1})^{v_{1}}\cdot (X_{n}-c_{2})^{v_{2}}(r_{1}\cdot pp1+r_{2}\cdot pp2)+the~ rest$

$\qquad~~~~=X_{0}\cdot (X_{n}-c_{1})^{v_{1}}\cdot (X_{n}-c_{2})^{v_{2}}+the~ rest$

None monomial in $the~rest$ is greater than or equal to $X_{0}$ , so the leading term of $f$ is apparently $X_{0}\cdot X_{n}^{v}$ which is equal to $ed$. Moreover we naturally have the following \textbf{Proposition 1} for $i=2$.

\textbf{Proposition 1:}
For every polynomial $f$ we get in the algorithm, $LC_n(f)=\prod_{j=1}^{i}(X_{n}-c_{j})^{v_{j}}$.

When $i>2$, assume the \textbf{Proposition 1} holds for $i-1$. $\forall~ ed\in E(D)$, $v:=proj(ed)$ and $X_0:=\frac{X^{ed}}{X_n^{v}}$. According to Section 4, we have $v=v_{1}+\ldots+v_{i}$. Based on the \textbf{Proposition 1}, \textbf{Remark 1} and \textbf{Remark 2}, $f_{1}$ and $f_{2}$ can be written as polynomials of $k(X_{n})[X_{1},\ldots,X_{n-1}]:$ $f_{1}=X_{0}\cdot \prod_{j=1}^{i-1}(X_{n}-c_{j})^{v_{j}}+the~rest$ and $f_{2}=X_{0}\cdot (X_{n}-c_{i})^{v_{i}}+the~rest$ and none of the monomials in $the~rest$ is greater than or equal to $X_{0}$.
Because $f_{1}$ and $\prod_{j=1}^{i-1}(X_{n}-c_{j})^{w_{j}+1}$ vanish on $\bigcup_{j=1}^{i-1}H_{j}$, $f_{2}$ and $(X_{n}-c_{i})^{w_{i}+1}$ vanish on $H_{i}$, we know that $q_{1}=f_{1}\cdot (X_{n}-c_{i})^{w_{i}+1}$ and $q_{2}=f_{2}\cdot \prod_{j=1}^{i-1}(X_{n}-c_{j})^{w_{j}+1}$ both vanish on $\bigcup_{j=1}^{i} H_{j}$. Then $f$ vanishes on $\bigcup_{j=1}^{i} H_{j}$ where $f=r_{1}\cdot q_{1}+r_{2}\cdot q_{2}$.

$$f=X_{0}\cdot\prod_{j=1}^{i}(X_{n}-c_{j})^{v_{j}}(r_{1}\cdot (X_{n}-c_{i})^{w_{i}+1-v_{i}}+r_{2}\cdot \prod_{j=1}^{i-1}(X_{n}-c_{j})^{w_{j}+1-v_{j}})+the~ rest$$

$\qquad\qquad=X_{0}\cdot\prod_{j=1}^{i}(X_{n}-c_{j})^{v_{j}}(r_{1}\cdot pp1+r_{2}\cdot pp2)+the~ rest$

$\qquad\qquad=X_{0}\cdot\prod_{j=1}^{i}(X_{n}-c_{j})^{v_{j}}+the~ rest$

None monomial in $the~rest$ is greater than or equal to $X_{0}$ and the leading term of $f$ is apparently $X_{0}\cdot X_{n}^{v}$ which is equal to $ed$. Hence the \textbf{Proposition 1} holds for arbitrary $i$.

Therefore we have proved that for any element $ed\in E(D)$, $f_{ed}:=f$ vanishes on $H$ and the leading term is $ed$. In the algorithm, we compute $f_{ed}$ in turn according to the lexicographic ordering of the elements of $E(D)$. Once we get a polynomial, we use the polynomials we got previously to reduce it ([step 4.6]). Now to prove the lower terms of $f^{'}$ are all in $D$ after such a reduction operation.

Let $D$ be a lower set, $a$ be a monomial, define $L(a,D)=\{b\in\mathbb{N}_{0}^{n};b\prec a,b\in D\}$. Given any $d\notin D$, there exist only two situations:
$d\in E(D)$ or $d\notin E(D)$ but $\exists d^{'}\in E(D),~ s.t.~d^{'}$ is a factor of $d$. Of course $d^{'}\prec d$.

The very first vanishing polynomial we got in the algorithm is an univariate polynomial of $X_{n}$ with leading term being $T$. It is easy to check it's lower terms are in $D$. Since the polynomial is a vanishing polynomial, we can say that $T$ can be represented as the linear combination of the elements of $L(T,D)$.

Since $T$ is the first element which is not in $D$ under lexicographic ordering. We \textbf{assume} that there exists such a monomial $M\notin D(M\succ T)$ that $\forall m\prec M(m\notin D)$, $m$ can be represented as the linear combination of the elements of $L(m,D)$. Now to prove $M$ could be represented as the linear combination of the elements of $L(M,D).$

If $M\in E(D)$, then the algorithm provides us a vanishing polynomial whose leading term is $M$ i.e. that $M$ can be represented as the combination of the terms which are all smaller than $M$. According to the assumption, for any lower term $m~(m\notin D)$ of the polynomial, $m$ can be represented as the linear combination of the elements of $L(m,D)$, then $M$ could be represented as the linear combination of the elements of $L(M,D).$

If $M\notin E(D)$, there exists $d^{'}\in E(D)$ s.t. $M=M^{'}\cdot d^{'}$. Since $d^{'}\prec M$, according to the assumption, we can substitute $d^{'}$ with the linear combination of the elements of $L(d^{'},D)$. Since all the elements in $L(d^{'},D)$ are smaller than $d^{'}$, then $M$ could be represented as the combination of elements which are all smaller than $M$. Then for the same reason described in the last paragraph, $M$ could be represented as the linear combination of the elements of $L(M,D).$

Therefore specially for any $ed\in E(D)$, all the lower terms of the polynomial $f_{ed}$ we got in the algorithm after the reduce operation are in $D$, and the proof is done.

\textbf{Theorem 4:} Given a set of points $H$ with multiplicity structures, $poly(H)$ is the reduced Gr$\ddot{\rm{o}}$bner basis of the vanishing ideal $I(H)$ and $D(H)$ is the quotient basis under lexicographic ordering.

\textbf{Proof:}
Let $m$ be the number of functionals defined by $H$ and then $m=dim(k[X]/I(H))$. Denote by $J$ the ideal generated by $poly(H)$. According to the \textbf{Assumption}, $poly(H)\subseteq I(H)$. So $dim(k[X]/I(H))\leq dim(k[X]/J)$. Let $C$ be the leading terms of polynomials in $J$ under lexicographic ordering, then $C\supseteq \bigcup_{\beta\in E(D(H))}(\beta+\mathbb{N}_{0}^{n})$ where the latter union is equal to $\mathbb{N}_{0}^{n}-D(H)$. Then we can get $C^{'}=\mathbb{N}_{0}^{n}-C\subseteq D(H)$. Because $k[X]/J$ is isomorphic as a $k$-vector space to the $k$-span of $C^{'}$, here $C^{'}$ is viewed as a monomial set. So $dim(k[X]/J)\leq \sharp D(H)=m$. Hence we have $$m=dim(k[X]/I(H))\leq dim(k[X]/J)\leq m.$$

Therefore $J=I(H)$, where $J=\langle poly(H)\rangle$. Hence apparently $poly(H)$ is exactly the reduced Gr$\ddot{\rm{o}}$bner basis of the vanishing ideal under lexicographic ordering, and $D(H)$ is the quotient basis.

\section{Intersection of ideals and some applications}

Some steps of our algorithm actually do the work of computing the intersection of two ideals, but we note that the information of the zeros of the ideals is necessary there (see [step 4.1] - [step 4.7] of \textbf{p-General case} in Section 6). We now bring up a new algorithm to compute the intersection of two ideals which does not require the information of the zeros of the ideals.

\textbf{Lemma 1:} $G$ is the reduced Gr$\ddot{\rm{o}}$bner basis of some $n$-variable polynomial ideal under lexicographic ordering with $X_{1}\succ X_{2}\succ\ldots\succ X_{n}$. Define $p_{0}(G)$ as the  univariate polynomial of $X_{n}$ in $G$.  View $g\in G$ as polynomial of $K(X_{n})[X_{1},\ldots,X_{n-1}]$ and define $LC_n(g)$ to be the leading coefficient of $g$ which is an univariate polynomial of $X_{n}$ and we have the conclusion that $LC_n(g)$ is always a factor of $p_{0}(G)$.

\textbf{Proof:} In fact \textbf{Proposition 1} in Section 6 holds for any given reduced Gr$\ddot{\rm{o}}$bner basis under lexicographic ordering since it is unique and can be constructed in the way our algorithm offers. According to the proposition, $\forall f\in G$, $LC_n(g)=\prod_{j=1}^{s}(X_{n}-c_{j})^{v_{j}}$ and $v_{j}\leq w_{j}+1.$ $p_{0}(G)=\prod_{j=1}^{s}(X_{n}-c_{j})^{w_{j}+1}$. Hence the proof is done.

Based on \textbf{Proposition 1} and \textbf{Lemma 1}, we give the algorithm \textbf{Intersection} to compute the intersection of two ideals $I_1$ and $I_2$ which are represented by the lexicographic ordering reduced Gr$\rm{\ddot{o}}$bner bases $G_{1}$ and $G_{2}$ and the greatest common divisor of $p_{0}(G_{1})$ and $p_{0}(G_{2})$ equals to 1. Denote by $Q(G)$ the quotient basis where $G$ is the reduced Gr$\rm{\ddot{o}}$bner basis. Algorithm \textbf{GP} is a sub-algorithm called in algorithm \textbf{Intersection}.

Algorithm \textbf{GP:} $G$ is a reduced Gr$\ddot{\rm{o}}$bner basis, for any given monomial $LT$ which is not in $Q(G)$, we get a polynomial $p$ in $\langle G\rangle$ whose leading term is a factor of $LT$: the $X_{1},\ldots,X_{n-1}$ components of the leading term are the same with that of $LT$ and the $X_{n}$ component has the lowest degree. Denoted by $p:=GP(LT,G).$

[step 1:] $G^{'}:=\{g\in G|$ the leading monomial of $g$ is a factor of $LT$ $\}$.

[step 2:] $G^{''}:=\{g\in G^{'}|\nexists g^{'}\in G^{'},~s.t.~$ the degree of $X_{n}$ of the leading monomial of $g^{'}$ is lower than that of $g$ $\}$.

[step 3:] Select one element of $G^{''}$ and multiply it by a monomial of $X_{1},\ldots,X_{n-1}$ to get $p$ whose leading monomial is $LT$.

Algorithm \textbf{Intersection:} $G_{1}$ and $G_{2}$ are the reduced Gr$\ddot{\rm{o}}$bner bases of two different ideals satisfying that $GCD(p_{0}(G_{1}),p_{0}(G_{2}))=1$. Return the reduced Gr$\ddot{\rm{o}}$bner basis of the intersection of these two ideals, denoted by $G:=Intersection(G_{1},G_{2})$.

[step 1:] $D:=Q(G_{1})+Q(G_{2})$. View $E(D)$ as a monomial set. $G:=\varnothing$.

[step 2:] If $E(D)=\varnothing$, the algorithm is done. Else select the minimal element of $E(D)$, denoted by $T$. $E(D):=E(D)/\{T\}$.

[step 3:] $$f_{1}:=GP(T,G_{1}),~f_{2}:=GP(T,G_{2}).$$$$q_{1}:=f_{1}\cdot p_{2},~q_{2}:=f_{2}\cdot p_{1}.$$

[step 4:] $$t_{1}:=\frac{p_{0}(G_2)}{LC_n(f_{2})},~t_{2}:=\frac{p_{0}(G_1)}{LC_n(f_{1})}.$$

[step 5:] Use Extended Euclidean Algorithm to find $r_{1},r_{2}$ s.t.
$$r_{1}\cdot t_{1}+r_{2}\cdot t_{2}=1.$$

[step 6:] $f:=q_{1}\cdot r_{1}+q_{2}\cdot r_{2}$. Reduce $f$ with $G$ to get $f^{'}$, and $G:=G\bigcup\{f^{'}\}$. Go to [Step 2].

Because the algorithm is essentially the same with [step 4.1] - [step 4.7] of \textbf{p-General case} in Section 6, here we don't give the proof.

The \textbf{Proposition 1} and \textbf{Lemma 1} reveal important property of the reduced Gr$\ddot{\rm{o}}$bner basis under lexicographic ordering. If a set of polynomials does not have this property, it is surely not a reduced Gr$\ddot{\rm{o}}$bner basis.

It is well-known that the Gr$\ddot{\rm{o}}$bner basis of an ideal under lexicographic ordering holds good algebraic structures and hence is convenient to use for polynomial system solving. To compute the zeros of an zero dimensional ideal with the reduced Gr$\ddot{\rm{o}}$bner basis $G$, we need first compute the roots of $p_0(G)$. Since $LC_n(g)$ ($g\neq p_0(G),~g\in G$) is a factor of $p_0(G)$, compute the roots of $LC_n(g)$ which has a smaller degree would be helpful for saving the computation cost.

\section{Conclusion}

Based on the algorithm \textbf{Intersection} in Section 7, the algorithm of \textbf{p-General case} in Section 6 can be simplified. The last sentence in [step 2] can be deleted and we can replace [step 4.3] and [step 4.4] by:

[step 4.3]:
$$f_{1}:=GLP(LT,D,G),f_{2}:=GLP(LT,D_{i}^{'},G_{i}).$$

[step 4.4]:
$$q_{1}:=f_{1}\cdot p_0(G_i);~~q_{2}:=f_{2}\cdot p_0(G).$$

$$pp1:=\frac{p_0(G_i)}{LC_n(f_{2})};~~pp2:=\frac{p_0(G)}{LC_n(f_{1})}.$$

During the induction of the algorithm in Section 6, we can record the leading coefficients for later use to save the computation cost and the computation cost is mainly on the Extended Euclid Algorithm. However the advantage of our algorithm is not fast computation, after all it depends on how many times we need to use the Extended Euclid Algorithm.

Our algorithm has an explicit geometric interpretation which reveals the essential connection between the relative position of the points with multiplicity structures and the quotient basis of the vanishing ideal. The algorithm offers us a new perspective of view to look into the reduced Gr$\ddot{\rm o}$bner basis which can help us understand the problem better. \textbf{Lemma 1} and the algorithm to compute the intersection of two ideals are the direct byproducts of our algorithm.

Since we finished the paper [1] previously which gives an algorithm to get the minimal monomial basis of Birkhoff interpolation problem with little computation cost, we have always believed that the algorithm could be interpreted in a more geometric way and the proof should be more beautiful and much easier to understand. The proof in [1] is so complicated that we ourselves don't like it. And it would be great if we can get the interpolation polynomial with little computation cost instead solving the linear equations since the minimal monomial basis can already be got in a simple way. That's why we began to study the vanishing ideal of the set of points with multiplicity structures which is essentially a special case of Birkhoff interpolation problem.

I still remember the moment when I first read the paper [2] written by Mathias Lederer in which the quotient basis and Gr$\ddot{\rm{o}}$bner basis can be got in a geometric way. I told myself that this was just what we wanted. Paper [2] concentrates on the the vanishing ideal of the set of points with no multiplicity structures in affine space. Although whether or not the points are with multiplicity structures matters much, the paper really inspired us a lot. Our algorithm also uses induction over variables and the definition of addition of lower sets is essentially the same with that in paper [2]. However during the induction procedure, we have to consider \textbf{p-Special case} and \textbf{p-General case}. This consideration, on one hand, clearly indicates the geometric meaning of the multiplicity structures of points, on the other hand, means a lot for our capacity of applying the algorithm of intersection of two ideals. In paper [2], the author uses Lagrange interpolation method to get the vanishing polynomial of all points from the polynomials vanishing on subsets of the points. However the Lagrange interpolation method could just not work for our problem because the points are with multiplicity structures. In this paper, we creatively use the Extended Euclidean Algorithm. Thanks goes to paper [2] and the author, after we solved the problem of the vanishing ideal of the set of points with multiplicity structures, we will move on to the Birkhoff problem.

\section{References}

[1]Na Lei, Junjie Chai, Peng Xia, Ying Li, A fast algorithm for multivariate Birkhoff interpolation problem, Journal of Computational and Applied Mathematics 236 (2011) 1656-1666.

[2]Mathias Lederer, The vanishing ideal of a finite set of closed points in affine space, Journal of Pure and Applied Algebra 212 (2008) 1116-1133.

[3]Hans J.Stetter. Numerical Polynomial Algebra. Chapter 2. SIAM, Philadelphia, PA, USA. 2004.

[4]M.G. Marinari, H.M. M$\rm{\ddot{o}}$ller, T. Mora, Gr$\rm{\ddot{o}}$bner bases of ideals defined by functionals with an application to ideals of projective points, J. AAECC 4 (2) (1993) 103-145.

[5]B$\acute{\rm{a}}$lint Felszeghy, Bal$\acute{\rm{a}}$zs R$\acute{\rm{a}}$th, Laj$\acute{\rm{o}}$s Ronyai, The lex game and some applications, Journal of Symbolic Computation 41 (2006) 663-681.

[6]L. Cerlinco, M. Mureddu, From algebraic sets to monomial linear bases by means of combinatorial algorithms, Discrete Math 139 (1995) 73-87.

\end{document}